\documentclass[12pt]{article}
\usepackage{amsmath,amsfonts,amsthm,enumerate}
\usepackage[utf8]{inputenc}
\usepackage{enumerate}
\usepackage{hyperref}
\newtheorem{theorem}{Theorem}
\newtheorem{lemma}[theorem]{Lemma}
\theoremstyle{remark}
\newtheorem*{remark}{Remark}

\newcommand{\Prob}{\mathbb{P}}
\newcommand{\Ex}{\mathbb{E}}

\author{Andreas N.\ Lager{\aa}s\footnote{supported by the Faculty of Science, G{\"o}teborg University.} \\  \small Mathematical Sciences and Centre for Theoretical Biology,\\  \small Chalmers University of Technology and G{\"o}teborg University, 412 96 Gothenburg, Sweden\\
\and Mathias Lindholm\footnote{supported by the Swedish Foundation for Strategic Research (SSF).} \\ \small Department of Mathematics, Stockholm University, 106 91 Stockholm, Sweden}
\date{\small Mathematics Subject Classification: 05C80}

\begin{document}

\title{A note on the component structure in random intersection graphs with tunable clustering}

\maketitle

\thispagestyle{empty}

\begin{abstract}
We study the component structure in random intersection graphs with tunable clustering, and show that the average degree works as a threshold for a phase transition for the size of the largest component. That is, if the expected degree is less than one, the size of the largest component is a.a.s.\ of logarithmic order, but if the average degree is greater than one, a.a.s.\ a single large component of linear order emerges, and the size of the second largest component is at most of logarithmic order.
\end{abstract}

\section{Introduction}\label{intro}

The random intersection graph, denoted $\mathcal{G}_{m, p}^{(n)}$, with a set of vertices $\mathcal{V} = \{v_{1}, \ldots, v_{n}\}$ and a set of edges $\mathcal{E}$ is constructed from a bipartite graph $\mathcal{B}^{(n)}_{m,p}$ with two sets of vertices: $\mathcal{V}$, identical to those of $\mathcal{G}_{m, p}^{(n)}$, and $\mathcal{A} = \{a_{1}, \ldots, a_{m}\}$, which we call auxiliary vertices. Edges in $\mathcal{B}^{(n)}_{m,p}$ between vertices and auxiliary vertices are included independently with probability $p \in [0, 1]$. An edge between two vertices $v_{i}$ and $v_{j}$ in $\mathcal{G}^{(n)}_{m,p}$ is only present in $\mathcal{E}$ if both $v_{i}$ and $v_{j}$ are adjacent to some auxiliary vertex $a_{k}$ in $\mathcal{B}^{(n)}_{m,p}$. Along the lines of Karo{\'n}ski et al. \cite{KSSC99} we set $m := \lfloor \beta n \rfloor$ and $p := \gamma n^{-(1+\alpha)/2}$, where $\alpha, \beta, \gamma \ge 0$, to obtain an interesting graph structure and bounded average vertex degree. For random (multi)graphs, the vertex degree distribution is defined as the distribution of the degree, i.e.\ the number of adjacent edges, of a vertex chosen uniformly at random. As has been shown by Stark \cite{S04}, the vertex degree distribution of the random intersection graph is highly dependent on the value of $\alpha$, but as shown by Deijfen and Kets \cite{DK07}, the clustering is {\it tunable} only when $\alpha \equiv 1$. In a recent paper by Behrisch \cite{B07}, the component structure of the random intersection graph is studied for $\alpha \neq 1$ and $\beta=1$, and the aim of the present note is to describe the component structure when $\alpha = 1$. We will henceforth keep $\beta$ and $\gamma$ fixed and positive, and sometimes suppress the dependency on these parameters in the notation: $\mathcal{G}^{(n)}$.

\section{The degree distribution}\label{egenskaper}

We define $D({m,n,p})$ to be a random variable with the vertex degree distribution of the graph $\mathcal{G}^{(n)}_{m,p}$. Stark \cite[Thm.\ 1]{S04} showed that the distribution of $D({m,n,p})$ has the following generating function
\begin{align*}
g_{D({m,n,p})}(z) := \Ex\left[z^{D({m,n,p})}\right] = \sum_{j = 0}^{n-1}\binom{n-1}{j}z^{j}(1-z)^{n-1-j}\left[1-p + p\left(1 - p\right)^{n-1-j}\right]^{m}.
\end{align*}
This distribution is from here onwards denoted ${\rm RIG}(m, n, p)$. Let us define a certain compound Poisson random variable $Z$ by its generating function
\begin{align*}
g_{Z}(s) := \Ex\left[s^{Z}\right] = \exp\big\{\lambda' \big(e^{\lambda''(s-1)} - 1\big)\big\},
\end{align*}
and write $Z\in{\rm CPoisson}(\lambda', \lambda'')$. Here $\Ex[Z]=\lambda'\lambda''$. Another result by Stark [6, Thm.\ 2], here slightly generalised, is
\begin{lemma}\label{graden}
If $n'$ and $n''$ are functions of $n$ such that $\lfloor \beta n\rfloor \geq n'$, $n'/n=\beta+o(1)$, $n\geq n''$, $n''/n=1+o(1)$, then $D({n',n'',\gamma/n})\overset{d}{\to}\mathrm{CPoisson}(\beta\gamma,\gamma)$ as $n\to\infty$.
\end{lemma}
This can be shown by inspecting the generating functions. In particular
\begin{align*}
\Ex[D({m,n,p})] &= g_{D({m,n,p})}'(1) = (n-1)[1-(1-p^{2})^{m}]\\
\Ex[D({m,n,p})(D({m,n,p})-1)] &= g_{D({m,n,p})}''(1)\\
&= (n-1)(n-2)[1 - 2(1-p^{2})^{m} + (1-p^{2}(2-p))^{m}],
\end{align*}
and with $n'$ and $n''$ as in the lemma we can deduce $\Ex[D(n',n'',\gamma/n)] = \mu + o(1) = O(1)$, where $\mu:=\beta\gamma^2$, and $\Ex[D(n',n'',\gamma/n)^{2}] = \mu(1+\mu+\gamma)+o(1)=O(1)$. We write 
\begin{align*}
g(s):= \exp\Big\{\beta\gamma\left(e^{\gamma(s-1)} - 1\right)\Big\}
\end{align*}
for the generating function of the limiting distribution $\mathrm{CPoisson}(\beta\gamma,\gamma)$.  Finally, let us define $\rho$ to be the smallest non-negative root of $\rho=g(\rho)$. 

\section{Results}\label{resultat}

\begin{theorem}\label{komponenten}
Let $\mu:=\beta\gamma^2$, i.e.\ the asymptotic expected degree of a randomly chosen vertex of $\mathcal{G}^{(n)}$.
\begin{enumerate}[(i)]
\item\label{sub} If $\mu < 1$, then there is a.a.s.\ no connected component in $\mathcal{G}^{(n)}$ with more than $O(\log n)$ vertices.
\item\label{sup} If $\mu > 1$, then $0<\rho<1$ and there exists a unique giant component of size $(1-\rho+o_{p}(1))n$, and the size of the second largest connected component is a.a.s.\ no larger than $O(\log n)$.
\end{enumerate}
\end{theorem}
With $W_n=o_p(a_n)$ we mean that $W_n/a_n\to0$ in probability as $n\to\infty$. As mentioned in the introduction, Behrisch has investigated the component structure for the random intersection graph when $\alpha \neq 1$, see \cite[Thm.\ 1]{B07}. It is worth noting that the results in Theorem \ref{komponenten} in this note are closer to the results that Behrisch obtained for the case $\alpha > 1$, than for $\alpha < 1$. For $\alpha < 1$ the size of the giant component is no longer linear in $n$.

\section{Proof of Theorem \ref{komponenten}}\label{beviset}

For the remainder of this note we will follow the notation and steps of the proof of Theorem 5.4 in \cite[Ch.\ 5.2]{JLR00}. Therefore most of the details that have not been altered from the original proof will be omitted. The proof is based on choosing a vertex at random from $\mathcal{V}$, say $v$, and exploring its component, say $\mathcal{C}(v)$. We start by \emph{visiting} the chosen vertex $v$ and \emph{identifying} its neighbours. Then we proceed by visiting an identified but unvisited vertex, if any remains, and identify its neighbours, and repeat this procedure until all vertices in the component have been visited. Let $X_{i}$ denote the number of newly identified vertices at the $i$th step of this exploration process. The event $\{|\mathcal{C}(v)|=k\}$ is equivalent to $\sum_{i=1}^kX_i=k-1$, and it is thus important to understand the growth of this partial sums process. The random variables $X_1,X_2,\dots$ are not i.i.d.\ but the partial sums process can nevertheless be related to other partial sums processes with i.i.d.\ summands so that we obtain bounds on events of the type above. We will need the following result for these partial sums processes.
\begin{lemma}\label{svansen}
Let $\delta>0$ and $\tilde{X}:=\tilde{X}_1+\cdots+\tilde{X}_k$, where $\tilde{X}_1,\dots$ are i.i.d. as $D(n',n'',\gamma/n)$ of Lemma \ref{graden}. Then, for large enough $n$, there exists a positive constant $C:=C(\beta,\gamma,\delta)$ such that $\Prob(\tilde{X} \geq (1+\delta)\mu k)\leq e^{-Ck}$ and $\Prob(\tilde{X} \leq (1-\delta)\mu k)\leq e^{-Ck}$.
\end{lemma}

\begin{remark}
This bound on the tail probabilities works since $\tilde{X}$ is a sum of $k$ independent random variables. As $n\to\infty$, the RIG-distribution of the summands does not change much: It is more or less $\mathrm{CPoisson}(\beta\gamma,\gamma)$, which is a ``well behaved'' distribution, and as $k$ increases, we expect an exponential decay of probabilities away from the mean of the sum. Since $C$ is not further specified, this bound is only useful as $k$ tends to infinity, which it may or may not do as a function of $n$.
\end{remark}

Before we prove the lemma, note that we can construct a \emph{multigraph} $\mathcal{H}^{(n)}_{m,p}$ from the same bipartite graph $\mathcal{B}_{m,p}^{(n)}$ as we used in the construction of $\mathcal{G}^{(n)}_{m,p}$, by letting the number of edges between $v_i$ and $v_j$ equal the number of auxiliary vertices $a_k$ that are adjacent to both $v_i$ and $v_j$. We denote with $\mathrm{RIMG}(m,n,p)$ the degree distribution of $\mathcal{H}^{(n)}_{m,p}$. $\mathrm{RIMG}(m,n,p)$ clearly dominates $\mathrm{RIG}(m,n,p)$, as we can obtain $\mathcal{G}^{(n)}_{m,p}$ from $\mathcal{H}^{(n)}_{m,p}$ by coalescing multiple edges between vertices into one single edge. $\mathrm{RIMG}(m,n,p)$ is a compound binomial distribution with generating function
$$
h(z) = (1-p+p(1-p+pz)^{n-1})^m,
$$
since, by construction, a vertex $v_i\in \mathcal{H}^{(n)}_{m,p}$ is connected to a $\mathrm{Binomial}(m,p)$ number of auxiliary vertices, each of which being connected to an independent $\mathrm{Binomial}(n-1,p)$ number of vertices in $\mathcal{V}\setminus \{v_i\}$.

The expected value of $\mathrm{RIMG}(\lfloor\beta n\rfloor,n,\gamma/n)$ is thus $\lfloor\beta n\rfloor(n-1)\gamma^2/n^2=\mu+O(1/n)=\Ex[D(\lfloor\beta n\rfloor,n,\gamma/n)]+O(1/n)$, so the expected difference in vertex degree between the multigraph and the ordinary graph is only $O(1/n)$. With $\eta^{(n)}$ defined as the difference in the total number of edges in $\mathcal{H}^{(n)}_{\lfloor\beta n\rfloor,\gamma/n}$ and $\mathcal{G}^{(n)}_{\lfloor\beta n\rfloor,\gamma/n}$, we have $\Ex[\eta^{(n)}]=O(1)$, by summing over all vertices. This will be used in the proof of Theorem \ref{komponenten}(\ref{sup}).

\begin{proof}[Proof of Lemma \ref{svansen}]
Note that $\Ex[e^{\theta \tilde{X}}]=\Ex[e^{\theta \tilde{X}_1+\cdots+\theta \tilde{X}_k}]=\Ex[e^{\theta \tilde{X}_1}]^k$. Let $s>0$. We have
\begin{align}
\Prob(\tilde{X}\leq (1-\delta)\mu k) &= \Prob(e^{-s\tilde{X}}\geq e^{-s(1-\delta)\mu k})\notag\\
&\leq e^{s(1-\delta)\mu k}\Ex[e^{-s\tilde{X}}]=\left(e^{s\mu-s\delta\mu}\Ex[e^{-s\tilde{X}_1}]\right)^k,\label{lower_tail}\\
\Prob(\tilde{X}\geq (1+\delta)\mu k) &= \Prob(e^{s\tilde{X}} \geq e^{s(1+\delta)\mu k})\notag\\
&\leq e^{-s(1+\delta)\mu k}\Ex[e^{s\tilde{X}}]=\left(e^{-s\mu-s\delta\mu}\Ex[e^{s\tilde{X}_1}]\right)^k,\label{upper_tail}
\end{align}
by Markov's inequality. Since $e^{-s\tilde{X}_1}\leq 1-s\tilde{X}_1+\frac12 s^2 \tilde{X}_1^2$ for $s>0$,
\begin{align*}
\Ex[e^{-s\tilde{X}_1}]&\leq 1-sE[\tilde{X}_1] + \tfrac12 s^2\Ex[\tilde{X}_1^2] = \exp\{\log(1-sE[\tilde{X}_1] + \tfrac12 s^2\Ex[\tilde{X}_1^2])\}\\
&=\exp\{-s\Ex[\tilde{X}_1]+O(s^2)\}=\exp\left\{-s\left(\mu + o(1) + O(s)\right)\right\}.
\end{align*}
The right hand side of \eqref{lower_tail} is thus $\exp\left\{-s\left(\delta\mu + o(1) + O(s)\right)k\right\}$, and we can fix a small $s$, such that for large enough $n$, $s(\delta\mu + o(1) + O(s))$ is positive (regardless of the value of $k$), and thus $\Prob(\tilde{X}\leq (1-\delta)\mu k) \leq  e^{-C'k}$ for some positive $C'$.

For the second part of the proof, let $\hat{X}\in \mathrm{RIMG}(n',n'',\gamma/n)$, so that $\hat{X}\geq_d \tilde{X}_1$.
\begin{align*}
\Ex[e^{s\tilde{X}_1}]&\leq \Ex[e^{s\hat{X}}] = \left(1-\tfrac{\gamma}{n}+\tfrac{\gamma}{n}\left(1-\tfrac{\gamma}{n}+\tfrac{\gamma}{n}e^s\right)^{n''-1}\right)^{n'} \\
&< \exp\left\{\gamma\tfrac{n'}{n}\left(e^{\gamma\frac{n''-1}{n}(e^s-1)}-1\right)\right\}\\
&= \exp\left\{\gamma(\beta+o(1))\left(e^{\gamma(1+o(1))s(1+O(s))}-1\right)\right\} \\
&= \exp\left\{\gamma\beta(1+o(1))\left(e^{\gamma s(1+o(1)+O(s))}-1\right)\right\} \\
&= \exp\left\{\mu s(1+o(1)+O(s))\right\}.
\end{align*}
The right hand side of \eqref{upper_tail} is thus less than $\exp\{-s(\delta\mu + o(1) + O(s))k\}$, and we can fix a small $s$, such that for large enough $n$, $s(\delta\mu + o(1) + O(s))$ is positive (regardless of the value of $k$), and thus $\Prob(\tilde{X}\geq (1+\delta)\mu k) \leq  e^{-C''k}$ for some positive $C''$. We conclude the proof of the lemma by letting $C=\min\{C',C''\}$.
\end{proof}

\begin{proof}[Proof of Theorem \ref{komponenten}(\ref{sub})] The process of exploring vertices that was briefly described in the beginning of Section \ref{beviset}, implies that $X_1$, the number of neighbours of the initially picked vertex, has distribution ${\rm RIG}(\lfloor \beta n\rfloor, n, \gamma/n)$. This, together with the fact that vertices only can be newly identified once, implies that $\sum_{i=1}^kX_i\leq_d \sum_{i=1}^kX^+_i$ for all $k$, where $X_1^+,\dots$ are i.i.d.\ ${\rm RIG}(\lfloor \beta n\rfloor, n, \gamma/n)$. Thus
\begin{align*}
\Prob(\exists i:|\mathcal{C}(v_i)|\geq k) &\leq \sum_{i = 1}^{n}\Prob\left(|\mathcal{C}(v_i)|\geq k \right) = n\Prob\left(|\mathcal{C}(v)|\geq k \right) \leq n\Prob\bigg( \sum_{j = 1}^{k}X_{j}^{+} \ge k-1 \bigg).
\end{align*}
Now we take $k:=k(n)$ increasing to infinity with $n$. Since all $X_{i}^{+}$ are i.i.d.\ \newline ${\rm RIG}(\lfloor \beta n\rfloor, n, \gamma/n)$, Lemma \ref{svansen} applies to $X^{+} := \sum_{j = 1}^{k}X_{j}^{+}$, which gives us
\begin{align*}
\Prob(\exists i:|\mathcal{C}(v_i)|\geq k)& \le n\Prob\left( X^{+} \ge k-1 \right) = n\Prob(X^{+} \ge (1+2\delta)\mu k - 1) \\
&\le n\Prob(X^{+} \ge (1+\delta)\mu k) \le n\exp\{-C k\},
\end{align*}
where $\mu < 1$, $\delta = (1/\mu - 1)/2 > 0$, $C$ is defined as in Lemma \ref{svansen}, and the penultimate inequality follows from $\delta\mu k>1$ for large enough $k$. That is, if $k(n) := \lceil (1+\epsilon) \log n / C\rceil$, $\epsilon > 0$, then $\Prob(\exists i:|\mathcal{C}(v_i)|\geq k) \le n^{-\epsilon} \to 0$ as $n \to \infty$, and the first part of Theorem \ref{komponenten} is proved.\end{proof}

\begin{proof}[Proof of Theorem \ref{komponenten}(\ref{sup})]

We will first show that there with probability tending to one are no clusters of size $k$ with $O(\log n) = k_-(n) \leq k \leq k_+(n):=n^{2/3}$. From now on, let $k_-\leq k \leq k_+$, where $k_-(n)$ will be specified shortly. The construction used in the proof is similar but more involved than the one of the first part of the proof. 

For the remainder of the proof we will implicitly condition on the event $\{\eta^{(n)}\leq \sqrt{n}\}$, whose probability tends to one when $n\to\infty$, by Markov's inequality and the fact that $\Ex[\eta^{(n)}]=O(1)$. Our construction fails on the complementary event, but this is of no consequence for the proof, since the probability of this event tends to zero. 

Let $A(v)$ be the event that the exploration process, initiated at $v$, at step $k_+$ has not terminated and that it at that step has identified fewer than $(\mu - 1)k_{+}/2$ vertices that have not yet been visited, i.e.\ $A(v)=\{k_+\leq \sum_{j=1}^{k_+}X_j\leq k_+-1 + (\mu - 1)k_{+}/2\}$. Let $B(v)$ be the event $\{\sum_{j=1}^{k_+}X_j\leq k_+-1 + (\mu - 1)k_{+}/2\}$. We will prove that the probability that the exploration process terminates after $k$ steps or that $A(v)$ holds for some $v$, tends to zero. Note that $\{|\mathcal{C}(v)|=k\}\subseteq B(v)$ for each $k$, and in particular $\{|\mathcal{C}(v)|=k_+\}\cup A(v)\subseteq B(v)$. We also have $\{|\mathcal{C}(v)|=k\}\subseteq \{\sum_{j=1}^{k}X_j\leq k-1 + (\mu - 1)k/2\}$ for each $k$.

On the set $B(v)\cap\{\eta^{(n)}\leq \sqrt{n}\}$, the exploration process has at step $k$ identified vertices in $\mathcal{V}$, that are adjacent to fewer than $(\mu+1)k_+/2+\sqrt{n}$ auxiliary vertices in $\mathcal{B}_{m,p}^{(n)}$. We claim that
\begin{align}\label{dominans}
\Prob\bigg( \sum_{i = 1}^{k}X_{i} \le k-1 + \frac{\mu - 1}{2}k \bigg) \leq \Prob\bigg( \sum_{i = 1}^{k}X_{i}^{-} \le k-1 + \frac{\mu - 1}{2}k \bigg)
\end{align}
holds with $X_1^{-},\dots,$ i.i.d.\ ${\rm RIG}(\lfloor \beta n - (\mu + 1)k_{+}/2-\sqrt{n}\rfloor, \lfloor n - (\mu + 1)k_{+}/2 \rfloor, \gamma/n)$. Note that $\sum_{i=1}^kX_i^-$  is \emph{not} a lower stochastic bound on $\sum_{i=1}^kX_i$ in the same way as $\sum_{i=1}^kX^+_i$ is an upper bound since the distribution of $X_1^-$ depends on $k_+$. The claim follows by a slight adaptation of the arguments of the proof of Theorem 4.3 in \cite[Ch.\ 4.2]{R08}: We compare our exploration process with another exploration process, which does not follow vertices that belong to a group of \emph{forbidden} vertices, or that are reached through edges generated by a group of forbidden auxiliary vertices. Both groups of forbidden vertices and auxiliary vertices are adjusted (diminished) after each step so that there are $(\mu+1)k_+/2$ vertices that are forbidden or identified, and $(\mu+1)k_+/2+\sqrt{n}$ auxiliary vertices that are forbidden or have generated an edge to an identified vertex. These adjustments are possible until the process has identified $(\mu+1)k_+/2$ vertices, which is long enough to deduce whether fewer than $k-1 + (\mu - 1)k/2$ vertices have been identified after step $k$. Furthermore, since we keep the number of forbidden vertices and auxiliary vertices fixed, the number of newly identified vertices by the modified exploration process will in each step be i.i.d.\ ${\rm RIG}(\lfloor \beta n - (\mu + 1)k_{+}/2-\sqrt{n}\rfloor, \lfloor n - (\mu + 1)k_{+}/2 \rfloor, \gamma/n)$.

Using \eqref{dominans}, assuming that $\{\eta^{(n)}\le\sqrt{n}\}$ holds, gives us that
\begin{align*}
\Prob(\exists i:\{k_-\leq |\mathcal{C}(v_i)|\leq k_+\}\cup A(v_i)) &\leq n\Prob(\{k_-\leq |\mathcal{C}(v)|\leq k_+\}\cup A(v)) \\
&= n\bigg(\sum_{k=k_-}^{k_+-1}\Prob(|\mathcal{C}(v)|=k)\!+\!\Prob(|\mathcal{C}(v)|=k_+)\!+\!\Prob(A(v))\bigg)\\ 
&\le n\sum_{k = k_{-}}^{k_{+}}\Prob\bigg(\sum_{j = 1}^{k}X_{j} \le k-1 + \frac{\mu - 1}{2}k\bigg)\\
&\le n\sum_{k = k_{-}}^{k_{+}}\Prob\bigg(\sum_{j = 1}^{k}X_{j}^{-} \le k-1 + \frac{\mu - 1}{2}k\bigg).
\end{align*}
We apply Lemma \ref{svansen} to $X^{-} := \sum_{j = 1}^{k}X_{j}^{-}$, which yields
\begin{align*}
\Prob(\exists i:\{k_-\leq |\mathcal{C}(v_i)|\leq k_+\}\cup A(v_i))&\leq n\sum_{k=k_-}^{k_+}\Prob\bigg(X^{-}\le k -1 + \frac{\mu-1}{2}k\bigg) \\
&= n\sum_{k=k_-}^{k_+}\Prob\bigg(X^{-}\le (1-\delta)\mu k-1\bigg)\\
&\leq n\sum_{k=k_-}^{k_+}\Prob(X^{-} \le (1-\delta)\mu k) \leq nk_+\exp\left\{-C k_-\right\},
\end{align*}
where $\mu > 1$, $\delta = (1-1/\mu)/2 > 0$ and $C$ is defined as in Lemma \ref{svansen}. Therefore, if $k_{-}(n) := \lceil(5/3 + \epsilon)\log n/C\rceil$, $\epsilon>0$ and $k_{+}(n) := \lfloor n^{2/3}\rfloor$ then $\Prob(\exists i:\{k_-\leq |\mathcal{C}(v_i)|\leq k_+\}\cup A(v_i))\leq n^{-\epsilon}\to 0$ as $n \to \infty$.

From Section \ref{intro} we know that two vertices in $\mathcal{G}^{(n)}$ are not connected if they avoid being adjacent to the same auxiliary vertex. Thus the probability that two vertices are not connected is $(1 - \gamma^{2}/n^{2})^{\lfloor \beta n\rfloor}$. Furthermore we know from the previous calculations that the probability that $A(v)$ holds for some $v$ tends to zero as $n$ tends to infinity, i.e.\ if there exist two different components of size $k_{+}$, they will each have at least $(\mu - 1)k_{+}/2$ identified but not yet visited vertices. This implies that the probability that two components each of size $k_{+}$ are disjoint after visiting their additional vertices is less than
\begin{align*}
\left( (1-\gamma^{2}/n^{2})^{\lfloor \beta n \rfloor} \right)^{((\mu - 1)k_{+}/2)^{2}} & \leq \exp\bigg\{-\frac{\gamma^{2}}{n^{2}} \lfloor\beta n\rfloor \bigg( \frac{\mu - 1}{2}\lfloor n^{2/3}\rfloor\bigg)^{2} \bigg\} \\
&\leq \exp\bigg\{-\frac{\mu (\mu - 1)^{2}}{4}O(n^{1/3})\bigg\}=o(1/n^2).
\end{align*}
That is, with probability tending to one, either vertices belong to connected components of size less than $k_{-}$, or to a \emph{unique} component of size at least $k_{+}$. 

To show that the size of the largest component grows linearly in $n$ with high probability, we need to show that the number of vertices that belong to small components, i.e.\ components of size $k_{-}$ or less, is strictly less than $n$, implying the remaining vertices belong to the giant component. Let $L_{i}:=\{|\mathcal{C}(v_i)|\leq k_-\}$, $Y_{i}:=\mathbf{1}_{L_i}$, and set $Y := \sum_{i = 1}^{n}Y_{i}$, so that $\Ex[Y]=n\Ex[Y_1]=n\Prob(L_1)$. By the same reasoning we use above, we can sandwich $\Prob(L_1)$ between $\Prob(C^+\leq k_-)$ and $\Prob(C^-\leq k_-)$ where $C^+$ and $C^-$ are the total sizes of branching processes with offspring distributed as $X^+_1$ and $X^-_1$, respectively. Lemma \ref{graden} implies that both offspring distributions tend to the same limit, $\mathrm{CPoisson}(\beta\gamma,\gamma)$, as $n$ tends to infinity. By standard results in branching process theory, see Athreya and Ney \cite[Thm.\ I.5.1]{AN72}, both probabilities $\Prob(C^+\leq k_-)$ and $\Prob(C^-\leq k_-)$ tend to the $\rho$ that we defined as the smallest non-negative root of $g(\rho)=\rho$, since $k_-(n)$ tends to infinity with $n$ and $\rho$ is the probability that the branching process with offspring distribution $\mathrm{CPoisson}(\beta\gamma,\gamma)$ has finite total size. It also holds that $0<\rho<1$, since $\mu>1$. Due to this, $\Ex[Y] = (\rho+o(1))n$, which implies that the expected size of the largest component is $(1-\rho + o(1))n$, and the proof that $Y$ is concentrated around $\rho n$ follows the last part of the proof of Theorem 1.(2) in Behrisch \cite[Sec.\ 4.2, p.\ 8]{B07} verbatim.\end{proof}

\section*{Acknowledgements}

We thank an anonymous referee for careful reading of the manuscript and for pointing out errors.

\end{document}